\documentclass[10pt]{amsart}

\usepackage{eucal}
\usepackage{amscd}
\usepackage[all,cmtip]{xy}
\usepackage{rotating}
\usepackage{hyperref,mathrsfs}
\usepackage{tikz}
\usepackage{amsmath,amssymb,soul}
\usepackage{multicol}
\usepackage[english,french]{babel}

\newtheorem{theorem}{Theorem }[section]

\newtheorem{lemma}[theorem]{Lemma}

\newtheorem{corollary}[theorem]{Corollary}
\newtheorem{proposition}[theorem]{Proposition}

\theoremstyle{definition}
\newtheorem{definition}[theorem]{Definition}
\newtheorem{remark}[theorem]{Remark}

\theoremstyle{abstract}

\def\eop{\hspace*{\fill}{\footnotesize$\blacksquare$}}

\newcommand{\wis}[1]{{\text{\em \usefont{OT1}{cmtt}{m}{n} #1}}}

\newcommand{\A}{\mathbb{A}}
\newcommand{\hP}{\mathbb{P}}
\newcommand{\B}{\overline{\mathbf{B}}}

\newcommand{\Fun}{\mathbb{F}_1}
\newcommand{\Spec}{\wis{Spec}}
\newcommand{\Z}{\mathbb{Z}}
\newcommand{\Proj}{\wis{Proj}}

\newcommand{\mF}{\mathcal{F}}

\newcommand{\fP}{\mathbf{P}}

\newcommand{\F}{\mathbb{F}}

\newcommand{\bA}{\mathbb{A}}
\newcommand{\bP}{\mathbb{P}}

\newcommand{\prf}{\textit{Proof. }}
\newcommand{\Sch}{\wis{Sch}}
\newcommand{\bL}{\mathbb{L}}
\usetikzlibrary{matrix}
\usetikzlibrary{arrows}
\usepackage{pgf}
\usepackage{lscape}
\usepackage{listings}

\usepackage{enumitem}

\title{Graphs, $\mathbb{F}_1$-schemes and virtual mixed Tate motives}
\author{Manuel M\'erida-Angulo and Koen Thas}
\address{Ghent University, Department of Mathematics, Krijgslaan 281, S22 and S25, B-9000 Ghent, Belgium}
\email{manmerang@gmail.com; koen.thas@gmail.com}

\begin{document}
\maketitle
\selectlanguage{english}
\begin{abstract}
In a number of recent works \cite{MMKT,MMKT1} the authors have introduced and studied a functor $\mF_k$ which associates to each loose graph $\Gamma$ |which is similar to a graph, but where edges with $0$ or $1$ vertex are allowed | a $k$-scheme, such that $\mF_k(\Gamma)$ is largely controlled by the combinatorics of $\Gamma$. Here, $k$ is a field, and we allow $k$ to be $\Fun$, the field with one element. For each finite prime field $\F_p$, it is noted in \cite{MMKT} that any $\mF_k(\Gamma)$ is polynomial-count, and the polynomial is independent of the choice of the field. In this note, we show that for each $k$, the class of $\mF_k(\Gamma)$ in the Grothendieck ring $K_0(\texttt{Sch}_k)$ is contained in $\Z[\bL]$, the integral subring generated by the virtual Lefschetz motive.
\end{abstract}
\selectlanguage{french}
\begin{abstract}
Dans certains ouvrages r\'{e}cents \cite{MMKT, MMKT1} les auteurs ont introduit et \'{e}tudi\'{e} un foncteur $\mF_k$ qui associe \`a chaque {\em loose graph} $\Gamma$ |similaire \`a un graph mais o\`u les ar\^etes avec 0 ou 1 sommets sont aussi permis | un $k$-sch\'ema, de telle fa\c{c}on que $\mF_k(\Gamma)$ est essentiellement controll\'e par la combinatoire de $\Gamma$. Ici, nous consid\'erons $k$ comme un corps, $\Fun$, le corps \`a un \'el\'ement, inclus. Pour chaque corps premier  $\mathbb{F}_p$, il est r\'emarqu\'e en \cite{MMKT} que tous les $\mF_k(\Gamma)$ sont {\em polynomial-count} et les polyn\^omes sont ind\'ependants du choix du corps. Dans cet article, on prouve que pour chaque $k$, la classe de $\mF_k(\Gamma)$ dans l'anneau de Grothedieck $K_0(\texttt{Sch}_k)$ est contenu dans $\Z[\bL]$, l'anneau int\'egral engendr\'e par le motif virtuel de Lefschetz.
\end{abstract}
\selectlanguage{english}




\section{Introduction}
\label{intro}

In a series of papers \cite{MMKT,MMKT1,MMKT2}, the authors of the present note have studied a functor $\mF_k$ which associates to each graph a $k$-scheme, where $k$ is any field. In fact, $\mF_k$ is a functor from the category of ``loose graphs'' | which are similar to graphs, but allowing edges with $0$ or $1$ vertex | and $k$ is any field, including the field ``with one element,'' $\mathbb{F}_1$. The functor $\mF_{\mathbb{F}_1}$ maps to Deitmar schemes, the core scheme theory over $\Fun$, where a slightly more general definition of Deitmar scheme is used than in \cite{Deitmarschemes2}, in order to give us more flexibility (see the definition of {\em congruence schemes} in \cite{Deitmarcongruence}). The functors obey a number of rules which model the fact that a projective $n$-dimensional $\Fun$-space  corresponds naturally to a complete graph on $n + 1$ vertices | see Tits \cite{anal} and Thas \cite{Chap1} | while an affine $n$-dimensional $\Fun$-space then logically corresponds to a loose graph on one vertex, with $n$ edges through it. The vertex corresponds to the closed point of the space, and the ``loose edges'' with affine directions | this idea is rooted in the formalism

\begin{equation}
\bP^n(k) = \bA^n(k) + \bP^{n - 1}(k). 
\end{equation}

The obtained schemes $\mF_k(\Gamma)$, with $\Gamma$ a loose graph, are covered by affine $k$-spaces whose intersections are governed by the relations in the loose graph. In \S \ref{Fk}, we will supply more details concerning the definition of $\mF_k$. 

One of the underlying ideas of \cite{MMKT,MMKT1} is that the schemes $\mF_k(\Gamma)$ can be studied by using the combinatorial theory of loose graphs, and that some (geometrical, topological) invariants can be easily determined as such. (This idea stems from the note \cite{KT-Japan}, where a related functor was defined, but lacked a for us important property that locally (loose) stars should always correspond to affine spaces.)
One example is that for each finite field $k = \mathbb{F}_q$, the number of $\mathbb{F}_q$-rational points of $\mF_k(\Gamma)$ can be easily determined if $\Gamma$ is a tree, solely from data of $\Gamma$. 
A second motivation is the fact that $\mF_k(\Gamma)$, for any loose graph $\Gamma$ and finite field $k$, is {\em defined over $\Fun$} in Kurokawa's sense \cite{Kurozeta}, and hence comes with a Kurokawa zeta function, as in \cite{Kurozeta}. This zeta function is independent of the choice of finite field,
and in this way, we introduce a new zeta function for all (loose) graphs which is very different than the Ihara zeta function.

\subsection{Polynomial countability and virtual mixed Tate motives}

In \cite{MMKT} it is shown that for each loose graph $\Gamma$, there exists a polynomial $P(X) \in \mathbb{Z}[X]$  such that for each finite field $k = \mathbb{F}_q \ne \mathbb{F}_1$
we have, with $\Big\vert \mF_{\F_q}(\Gamma) \Big\vert_q$ the number of $\F_q$-rational points of $\mF_{\F_q}(\Gamma)$, that
\begin{equation}
\Big\vert \mF_{\F_q}(\Gamma) \Big\vert_{q} = P(q).
\end{equation}

By definition, this means that the scheme $\mF_p(\Gamma)$  is {\em polynomial-count} for all primes $p$ \cite[Appendix]{Katz}, and moreover, the polynomial $P(X)$ is independent of the chosen finite field. 
It has been conjectured (and it would follow from one of the Tate conjectures) that this property implies that the class of $\mF_{\F_q}(\Gamma)$ in the Grothendieck ring  of schemes $K_0(\texttt{Sch}_{\F_q})$ is contained in $\Z[\bL]$, where $\bL := [\bA^1(\F_q)]$ is the class of the affine line, that is, that $[\mF_{\F_q}(\Gamma)]$ is a ``virtual mixed Tate motive.'' 

This question is the objective of the present note.

\subsection{The present note}

In this note, we obtain the following result.

\begin{theorem}
\label{main}
Let $\Gamma$ be any loose graph, and let $k \ne \Fun$ be any finite field. Then the class $[\mF_k(\Gamma)] \in K_0(\texttt{Sch}_k)$ is 
a virtual mixed Tate motive.
\end{theorem}

Our proof uses the process of ``surgery,'' which is a stepwise procedure devised in \cite{MMKT} and performed on a loose graph, in which each step consists of replacing an edge with $2$ vertices from a prescribed set of edges, by two edges with only one vertex. The local dimension of the corresponding ($k$-)scheme rises, and at the end of the process one winds up with a tree. By doing a local analysis, which refines a part of the ``Affection Principle'' obtained in \cite{MMKT}, we will show that the classes of $k$-schemes that arise in consecutive steps, always differ by elements in $\Z[\bL]$.
Finally, using precise results for trees, which were already obtained in \cite{MMKT}, we then finish the proof.

\section{Deitmar schemes}
\label{Deitmarintro}

We consider an ``$\Fun$-ring'' $A$ to be a multiplicative commutative monoid with an extra absorbing element 0.
Let $\Spec(A)$ be the set of all {\em prime ideals} of $A$ together with a Zariski topology. This topolo\-gical space endowed with a structure sheaf of $\Fun$-rings is called an {\em affine Deitmar scheme}, and is also denoted by $\Spec(A)$. We define a {\em monoidal space} to be a pair $(X, \mathcal{O}_X)$ where $X$ is a topological space and $\mathcal{O}_X$ is a sheaf of $\Fun$-rings defined over $X$. A {\em Deitmar scheme} is then a monoidal space such that for every point $x\in X$ there exists an open subset $U \subseteq X$ such that $(U, \mathcal{O}_X|_{U})$ is isomorphic to an affine Deitmar scheme.
For a more detailed definition of Deitmar schemes and the structure sheaf of $\Fun$-rings, we refer to \cite{Deitmarschemes2}.

\subsection{Affine space} 

Let $A:=\Fun[X_1, \ldots, X_n]$ be the {\em monoidal ring} in $n$ variables; then the {\em $n$-dimensional affine space over $\Fun$} is defined as the monoidal space $\Spec(A)$ and denoted by $\mathbb{A}^n_{\Fun}$ or $\bA^n(\Fun)$. There is one closed point, and to each unknown $X_i$
corresponds a coordinate axis (``direction'').

\subsection{Projective space} 

The $n$-dimensional projective space is defined as the projective scheme $\Proj(\Fun[x_0\ldots,x_n]) =: \bP^n_{\Fun}$, or $\bP^n(\Fun)$. There are $n + 1$ closed points, and a linear subspace of dimension $r$ contains $r + 1$ of these points; in particular, a projective subline has $2$ closed points. Combinatorially, one depicts $\bP^n(\Fun)$ as a complete graph on $n + 1$ vertices. The combinatorial affine $n$-space over $\Fun$ then arises by deleting a complete subgraph on $n$ vertices. So one obtains a loose graph with one vertex and $n$ loose edges.

\subsection{Embedding theorem} 
\label{EmbThm}

 Let $\Gamma$ be a loose graph. The embedding theorem of \cite{KT-Japan} observes that $\Gamma$ can be seen as a 
 subgeometry of the combinatorial projective $\Fun$-space $\mathbb{P}(\overline{\Gamma})$, where $\overline{\Gamma}$ and $\bP(\overline{\Gamma})$ are constructed as in section $\ref{Fk}$. Applying $\mF_k$ (for any field $k$ including $\Fun$), cf. the next section, one obtains that $\mF_k(\Gamma)$ is embedded in $\bP(\overline{\Gamma})$ (now seen as a scheme).

\section{The functors $\mF_k$}
\label{Fk}

We briefly describe how one can associate a Deitmar scheme to a loose graph $\Gamma$ through the functor $\mathcal{F}_k$, with $k$ any field, including $\Fun$.  The main thing is that the functor must obey a set of rules, namely: \medskip

\begin{itemize}
\item[COV]
If $\Gamma \subset \widetilde{\Gamma}$ is a strict inclusion of loose graphs, $\mF_k(\Gamma)$ also is a proper subscheme
of $\mF_k(\widetilde{\Gamma})$.
\item[LOC-DIM]
If $x$ is a vertex of degree $m \in \mathbb{N}^\times$ in $\Gamma$, then there is a neighborhood $\Omega$ of $x$ in 
$\mF_k(\Gamma)$ such that $\mF_k(\Gamma)_{\vert \Omega}$ is an affine space of dimension $m$.
\item[CO]
If $K_m$ is a sub complete graph on $m$ vertices in $\Gamma$, then $\mF_k(K_m)$ is a closed sub projective space 
of dimension $m - 1$ in $\mF_k(\Gamma)$. 
\item[MG]
An edge without vertices should correspond to a multiplicative group.
\end{itemize}
\medskip

Rule (MG) implies that we have to work with a more general version of Deitmar schemes since the multiplicative group $\mathbb{G}_m$ over $\Fun$ is defined  to be isomorphic to
\begin{equation*} \Spec(\Fun[X, Y]/(XY=1)),\end{equation*} 
where the last equation generates a congruence on the free abelian monoid $\Fun[X,Y]$. (The equation $XY = 1$ is not defined in Deitmar scheme theory.) 
The reader can find a more detailed explanation of this association in \cite{MMKT1}.

A very simple way to explain $\mF_k$ is as follows: first, for any loose star $S$ (one vertex plus a number $n$ of loose edges through it), $\mF_k(S)$ is the affine $k$-space of dimension $n$. Here, a {\em loose edge} is an edge with $0$ or $1$ vertex. 
Now let $\Gamma$ be any connected loose graph, and let $\overline{\Gamma}$ be the graph-theoretical completion of $\Gamma$, that is, $\overline{\Gamma}$ is the graph one obtains by adding a new vertex on any loose edge of $\Gamma$. Say that $\overline{\Gamma}$ has $m + 1$ vertices. Let $\bP(\overline{\Gamma})$ be the projective $\Fun$-space of dimension $m$ defined by these vertices (see \S \ref{Deitmarintro}); then $\mF_k(\Gamma)$ is the union in 
$\bP(\overline{\Gamma}) \otimes_{\Spec(\Fun)}\Spec(k)$ of the $k$-affine spaces defined
by the stars which are defined by each vertex, without the closed points which correspond to the vertices which were added to obtain $\overline{\Gamma}$. 
If we choose coordinates (in $\bP(\overline{\Gamma})$) such that each such vertex has as coordinates a vector in $\{ 0,1\}^{m + 1}$ with precisely one nonzero entry, then for each $k$, $\mF_k(\Gamma)$ can be described explicitly analytically. The disconnected case is easily derived from the connected case.

\begin{theorem}[\cite{MMKT1}]
The map $\mathcal{F}_{\Fun}$ is a functor from the category of loose graphs to the cate\-gory of Deitmar congruence schemes. Moreover, for any finite field $k$ (or $\mathbb{Z}$), the lifting map $\mathcal{F}_k(\cdot)=\mathcal{F}_{\Fun}(\cdot)\otimes_{\Spec(\Fun)} \Spec(k)$ is also a functor.  
\end{theorem}

Let $\Gamma$ be a loose graph and $\mathcal{F}_{\Fun}(\Gamma)$ be the Deitmar scheme associated to it. 
Let us call $v_1, \ldots, v_k$ the vertices of $\Gamma$ and $\A_{v_i}$ the affine space associated to $v_i$, $1\leq i\leq k$.

\begin{lemma}[\cite{MMKT}]
\label{l1}
For all $1\leq r,s, \leq k,$ $\A_{v_r}\cap\A_{v_s}\neq \emptyset$ if and only if $v_r$ and $v_s$ are adjacent vertices in $\Gamma$.
\end{lemma}

\section{Grothendieck ring of schemes of finite type over $\Fun$}

The $\Spec$-construction on monoids allow\-s us to have a scheme theory over $\Fun$ defined in an analogous way to the classical scheme theo\-ry over $\Z$. This also allows us to define the {\em Grothendieck ring of schemes over $\Fun$}.

\begin{definition}
The {\em Grothendieck ring of schemes of finite type over $\Fun$}, denoted as $K_0(\Sch_{\Fun})$, is generated by the isomorphism classes of schemes ${X}$ of finite type over $\Fun$, $[X]_{_{\Fun}}$, with the relation
\begin{equation}
[X]_{_{\Fun}}= [X\setminus Y]_{_{\Fun}} + [Y]_{_{\Fun}} 
\end{equation}
for any closed subscheme $Y$ of $X$ and with the product structure given by
\begin{equation}
[X]_{_{\Fun}}\cdot[Y]_{_{\Fun}}= [X\times_{\Fun}Y]_{_{\Fun}}.
\end{equation}
\end{definition}
 
We denote by $\underline{\bL}=[\mathbb{A}^1_{\Fun}]_{_{\Fun}}$ the class of the affine line over $\Fun$. \medskip

\section{Counting polynomials}

Let $\Gamma$ be a loose tree and $\mathcal{F}_{\Fun}(\Gamma)$ its corresponding Deitmar scheme.  The next result gives us information about the class of $\mathcal{F}_{\Fun}(\Gamma)$ in the Grothendieck ring of Deitmar schemes of finite type, $K_0(\Sch_{\Fun})$. 
We will sometimes use the notation $[\Gamma]_{_{\Fun}}$ for the class of $\mathcal{F}(\Gamma)$ in $K_0(\Sch_{\Fun})$ (also when $\Gamma$ is a general loose graph). We adapt the same notation over fields $k$.
\medskip

\begin{theorem}[\cite{MMKT}]
\label{D3.1}
Let $\Gamma$ be a loose tree. Let $D$ be the set of degrees $\{d_1, \ldots, d_k \}$ of $V(\Gamma)$ such that $1 < d_1 < d_2 < \ldots < d_k$ and let $n_i$ be the number of vertices of $\Gamma$ with degree $d_i$, $1\leq i \leq k$. We call $E$ the number of vertices of $\Gamma$ with degree 1 and $I= \sum_{i=1}^k n_i - 1$. Then 
\begin{equation}
\big[\Gamma\big]_{_{\Fun}} =  \displaystyle\sum_{i = 1}^k n_i\underline{\bL}^{d_i} - I\cdot\underline{\bL} + I + E.
\end{equation}
\end{theorem}


\section{Surgery} 

In order to inductively calculate the counting polynomial of a $\mathbb{Z}$-scheme
coming from a general loose graph, we introduced a procedure called {\em surgery}, in \cite{MMKT}. In each step of the procedure we ``resolve'' an edge, so as to eventually end up with a tree in much higher dimension. One has to keep track of how the counting polynomial changes in each step.

\subsection{Resolution of edges}

Let $\Gamma$ be a loose graph, and let $e$ be an edge with two distinct vertices $v_1, v_2$. The {\em resolution} of $\Gamma$
{\em along} $e$, denoted $\Gamma_e$, is the loose graph obtained from $\Gamma$ by deleting $e$, and adding two new
loose edges (each with one vertex) $e_1$ and $e_2$, where $v_i \in e_i$, $i = 1,2$.\medskip


The following theorem reduces the computation of the alteration of the number of $k$-rational points after resolving an edge, to a local problem.
In the statement, $\mathrm{d}(\cdot,\cdot)$ stands for the distance function in a graph.

\begin{theorem}[Affection Principle, \cite{MMKT}]
\label{AP}
Let $\Gamma$ be a finite connected loose graph, let $xy$ be an edge on the vertices $x$ and $y$, and let $S$ be a subset of the vertex set. 
Let $k$ be any finite field, and consider the $k$-scheme $\mF_k(\Gamma)$. Then $\cap_{s \in S}\A_s$,where $\A_s$ is the local affine space corresponding to the vertex $s\in S$, changes when one resolves 
the edge $xy$  only if $\cap_{s \in S}\A_s$ is contained in $\fP_{x,y}$, the projective subspace of $\mathbb{P}(\Gamma) \otimes_{\Fun}k$ generated by $\B(x,1) \cup \B(y,1)$, where $\B(x,1)=\{v\in V(\Gamma) ~|~ \mathrm{d}(v,x)\leq1\}$.
\end{theorem}

In terms of counting polynomials, we have the following theorem, in which $\Big\vert \cdot \Big\vert_k$ denotes the number of $k$-rational points. 

\begin{corollary}[Polynomial Affection Principle, \cite{MMKT}]
\label{PAP}
Let $\Gamma$ be a finite connected loose graph, let $\Gamma_{xy}$ be the loose graph after resolving the edge $xy$ and let $k$ be any finite field.
Then in $K_0(\Sch_k)$ we have
\begin{equation}\label{pap}
\Big\vert \Gamma \Big\vert_k - \Big\vert \Gamma_{xy} \Big\vert_k = \Big\vert \Gamma_{\vert \fP_{x,y}} \Big\vert_k - \Big\vert {\Gamma_{xy}}_{\vert \fP_{x,y}} \Big\vert_k.
\end{equation}
\end{corollary}

\medskip
\subsection{Counting polynomial for general loose graphs}\label{cplg}

To compute the counting polynomial of a scheme coming from a loose graph $\Gamma$ we choose a spanning loose tree $T$ of $\Gamma$ and resolve in $\Gamma$ all edges not belonging to $T$. This yields a loose tree $\overline{T}$ in which we apply the map defined in Theorem \ref{D3.1} so as to obtain a counting polynomial for it. Take an edge $e$ now that was resolved and consider the loose graph $\overline{T}^{e}$ in which all other edges except $e$ are resolved, i.e., $\overline{T}^{e}$ is the next-to-last step in the procedure of obtaining $\overline{T}$. Thanks to Corollary \ref{PAP}, we can compute the counting polynomial for $\overline{T}^{e}$ by restricting to $\fP_{e}$ (for the concrete formulas of the Affection Principle we refer to \cite[section 11]{MMKT}). By repeating this process as many times as edges were resolved, we inductively obtain the counting polynomial for $\mF_k(\Gamma)$. 

\begin{proposition}[\cite{MMKT}] 
Let $\Gamma$ be a loose graph and let $T$ and $\overline{T}$ be defined as above. Then the counting polynomial in $K_0(\Sch{\Fun})$ of $\mathcal{F}(\overline{T})$ is independent of the choice of the spanning loose tree $T$ of $\Gamma$.
\end{proposition}

\newpage
\section{The main result}

Let $\Gamma$ be any loose graph, and let $\mF_k(\Gamma)$ be the corresponding $k$-scheme, with $k$ any field (including $\mathbb{F}_1$).
Our proof goes by induction on the number $N$ which is the sum of the number of edges and the number of vertices. For small values of $N$, the main theorem, which is stated below for the sake of convenience, is easy to obtain.

\begin{theorem}
\label{main2}
Let $\Gamma$ be any loose graph, and let $k \ne \Fun$ be any finite field. Then the class $[\mF_k(\Gamma)] \in K_0(\texttt{Sch}_k)$ is 
a virtual mixed Tate motive.
\end{theorem}

We now proceed with the proof. Note that we may suppose w.l.o.g. that $\Gamma$ is connected (and note that resolving an edge on a connected loose graph not necessarily yields again a connected loose graph). Note also that resolution of edges on trees is not defined.

\bigskip
\subsection{Local lemmas}

The following lemma is easy to prove:

\begin{lemma}
\label{under}
In $K_0(\Sch_k)$, we have that $[\mF_k(\Gamma)] \in \Z[\bL]$ if and only if $[\mF_k(\underline{\Gamma})] \in \Z[\bL]$, with
$\underline{\Gamma}$ the underlying graph of $\Gamma$. \eop
\end{lemma}

The underlying graph of $\Gamma$ is the subgraph of $\Gamma$ induced on the vertices.
By Lemma \ref{under}, we may thus suppose that $\Gamma$ is a graph. Now suppose $e = xy$ is an edge, with $x$ and $y$ its vertices. 
Resolve the edge $xy$ to obtain $\Gamma_{xy}$ (this is a loose graph). \medskip

\begin{remark} Notice that intersecting with a projective space commutes with the functor $\mathcal{F}_k(\cdot)$. We will prove this remark in the following lemma.\end{remark}

\begin{lemma} Let us denote by $\fP=K_V$ the complete graph defined on a subset $V$ of vertices of $\Gamma$ and let us call $\mathbb{P}_k$ the $k$-projective space defined by $\fP$. Then $\mathcal{F}_k(\Gamma)\cap\mathbb{P}_k = \mathcal{F}_k(\Gamma\cap\fP)$.
\end{lemma}

\prf Before starting the proof let us denote by $S_w$ the {\em loose star} of a vertex $w$ of $\Gamma$, that is, the loose subgraph of $\Gamma$ formed by the vertex $w$ and all its incident edges.\medskip

It is easy to check that $\mathcal{F}_k(\Gamma\cap\fP)$ is a subscheme of $\mathcal{F}_k(\Gamma)\cap\mathbb{P}_k$ since $\Gamma\cap\fP$ is a subgraph of both $\Gamma$ and $\fP$. Consider now a point $x \in \mathcal{F}_k(\Gamma)\cap\mathbb{P}_k$. Then, from the definition of $\mF_k$ (see \S \ref{Fk}),  $x$ belongs to $\Spec(\A_v)\cap\mathbb{P}_k$, for a vertex $v\in \Gamma$. The latter scheme is defined by the part of the loose star $S_v \subseteq \Gamma$ inside $\fP$, i.e., by $S_v\cap\fP$. This concludes the proof since $S_v\cap\fP$ is a subgraph of $\Gamma\cap\fP$ and so $x\in \mF_k(S_v\cap\fP) \subseteq \mF_k(\Gamma\cap\fP)$.
\eop\medskip

The following lemma, in the spirit of Corollary \ref{PAP}, shows that we can restrict ourselves to local
considerations.

\begin{lemma}
\label{local}
In $K_0(\Sch_k)$, we have that $[\mF_k(\Gamma)] - [\mF_k(\Gamma_{xy})] \in \Z[\bL]$ if and only if
$[\mF_k(\Gamma \cap \fP_{x,y})] - [\mF_k(\Gamma_{xy} \cap \fP_{x,y})] \in \Z[\bL]$.
\end{lemma}

\prf In order to prove the statement, we will compute the difference of classes $[\mathcal{F}_k(\Gamma)] - [\mathcal{F}_k(\Gamma_{xy})]$. Thanks to the remark above and the relative topology on $\mathcal{F}_k(\Gamma)$ and $\mathcal{F}_k(\Gamma_{xy})$, we can deduce that both $\mathcal{F}_k(\Gamma\cap\fP_{x,y})$ and $\mathcal{F}_k(\Gamma_{xy}\cap\fP_{x,y})$ are closed in $\mathcal{F}_k(\Gamma)$ and $\mathcal{F}_k(\Gamma_{xy})$, respectively. Then, by the relations in the appropriate Grothendieck ring of schemes, we have that:
\begin{equation}\label{CLASS}\begin{cases}
[\mathcal{F}_k(\Gamma)] = [\mathcal{F}_k(\Gamma\cap\fP_{x,y})] + [\mathcal{F}_k(\Gamma)\setminus \mathcal{F}(\Gamma\cap\fP_{x,y})],\\
[\mathcal{F}_k(\Gamma_{xy})] = [\mathcal{F}_k(\Gamma_{xy}\cap\fP_{x,y})] + [\mathcal{F}_k(\Gamma_{xy}) \setminus \mathcal{F}(\Gamma_{xy}\cap\fP_{x,y})].
\end{cases}\end{equation}

We will prove that the last terms on the right-hand side of the equations are the same. Let $\Gamma'=\Gamma\cap\fP_{x,y}$ and $\Gamma'_{xy}=\Gamma_{xy}\cap\fP_{x,y}$. Note that thanks to the Affection Principle (see Theorem \ref{AP}, and \cite[Lemma 11.5]{MMKT}), in order to compare the classes of $\mF_k(\Gamma)\setminus\mF_k(\Gamma')$ and $\mF_k(\Gamma_{xy})\setminus\mF_k(\Gamma'_{xy})$ in $K_0(\Sch_k)$, we (only) need to take into account the local affine spaces in $\mF_k(\Gamma)$ ($\mF_k(\Gamma_{xy})$) associated to vertices of $\Gamma$ ($\Gamma_{xy}$) which are at distance at most one from the loose graph $\Gamma\setminus\Gamma'$ ($\Gamma_{xy}\setminus\Gamma'_{xy}$) (since vertices at distance strictly more than one, give rise to affine spaces that remain unchanged through resolution of $xy$). \medskip

From the definition of $\Gamma'$, it is easy to see that both vertices $x$ and $y$ are at least at distance two from any vertex of $\Gamma\setminus\Gamma'$, which implies that 
\begin{equation}
\label{eqlocal}
\begin{cases}
(\mathcal{F}_k(\Gamma)\setminus \mathcal{F}_k(\Gamma'))\cap\mathbb{A}_x=\emptyset,\\
(\mathcal{F}_k(\Gamma)\setminus \mathcal{F}_k(\Gamma'))\cap\mathbb{A}_y=\emptyset.\end{cases}\end{equation}

As resolving the edge $xy$ only changes locally the affine spaces $\A_x$ and $\A_y$ in $\mF(\Gamma)$ (more precisely in $\mF(\Gamma')$), and as the distance between $x$ (or $y$) and $\Gamma\setminus\Gamma'$ is preserved through resolution, this process does not affect the local affine spaces in $\mathcal{F}_k(\Gamma)\setminus \mathcal{F}_k(\Gamma')$, nor the intersection of any two of them. Notice that in the case of vertices $v\in\Gamma'$ at distance one from $\Gamma\setminus\Gamma'$, possible changes of the affine space $\A_v$ in $\mF(\Gamma)$ by resolution of $xy$ do not affect the scheme $\mathcal{F}_k(\Gamma)\setminus \mathcal{F}_k(\Gamma')$; changes only occur in the completion $\overline{\A_v}\cap\mF(\Gamma')$.\medskip

It is now easy to observe that there is a natural isomorphism  between $\mathcal{F}_k(\Gamma)\setminus \mathcal{F}_k(\Gamma')$ and $\mathcal{F}_k(\Gamma_{xy})\setminus \mathcal{F}_k(\Gamma'_{xy})$ induced by the graph morphism
\begin{equation*}
\gamma: \Gamma'' \rightarrow \Gamma_{xy}'',
\end{equation*} 

\noindent where $\Gamma''$ (respectively $\Gamma_{xy}''$) is the subgraph of $\Gamma$ (respectively $\Gamma_{xy}$) defined on the vertex set $V(\Gamma\setminus\Gamma')\cup\{v\in\Gamma' ~|~ d(v, \Gamma\setminus\Gamma')=1\}$ (respectively $V(\Gamma_{xy}\setminus\Gamma_{xy}')\cup\{v\in\Gamma_{xy}' ~|~ d(v, \Gamma_{xy}\setminus\Gamma_{xy}')=1\}$), and where $\gamma$ acts as the identity on vertices. This implies that both classes in $K_0(\Sch_k)$ are equal. We can conclude now that 
\begin{equation*}
[\mathcal{F}_k(\Gamma)] - [\mathcal{F}_k(\Gamma_{xy})] = [\mathcal{F}_k(\Gamma\cap\fP_{x,y})] - [\mathcal{F}_k(\Gamma_{xy}\cap\fP_{x,y})]. \end{equation*}\eop \\

By Lemma \ref{local}, we may suppose that $\Gamma = \Gamma \cap \fP_{x,y}$. 

We now slightly refine the Affection Principle from \cite{MMKT}.

\begin{lemma} 
\label{affect2}
Let $\Gamma$ be a graph, $xy$ an edge with vertices $x$ and $y$ and $\Gamma_{xy}$ the graph after resolving the edge $xy$. Let $u,v$ be two vertices of $\Gamma$ and consider $\A_u$ and $\A_v$, the local affine spaces at $u$ and $v$ in $\mF_k(\Gamma)$. The intersection  $\A_u\cap\A_v$ and the union $\A_u \cup \A_v$ change after resolution only if $u,v\in\{x,y\}\cup(x^{\perp}\cap y^{\perp})$.
\end{lemma}

\prf
First let us note that if $\A_u\cap\A_v=\emptyset$, then $\A_u\cup\A_v$ is stable under resolution. Consider now a vertex $u\in x^{\perp}\setminus ((x^{\perp}\cap y^{\perp})\cup\{y\})$. Then, it is clear that neither $\A_u$ nor $\overline{\A_u}\cap\mF_k(\Gamma)$ changes after resolving the edge $xy$. The latter is not affected by the resolution since the edge $xy$ is not in the graph $\Gamma\cap\B(u,1)$. The same holds for vertices $u\in y^{\perp}\setminus ((x^{\perp}\cap y^{\perp})\cup\{x\})$. To simplify notation we will write from now on only $\overline{\A_u}$ instead of $\overline{\A_u}\cap\mF_k(\Gamma)$ and we consider it embedded in the ambient space of $\mF_k(\Gamma)$.\medskip

Now suppose that $u\in x^{\perp}\cap y^{\perp}$; then the graph defined by $xy$ is a subgraph of the ``part at infinity'' of the graph completion of $S_u$, the star associated to $u$. This implies that locally at $u$ the changes that occur by resolving $xy$ are contained in $\overline{\A_u}\setminus\A_u$, so the local affine space at $u$ also remains invariant under resolution of $xy$.\medskip

Observe that $\A_u\cap\A_v$ and $\A_u\cup\A_v$ are controlled by $\A_u, ~\A_v, ~\overline{\A_u}, ~\overline{\A_v}$ and $\overline{\A_u}\cap\overline{\A_v}$. So, if $u,v\neq x,y$ and $u,v\notin x^{\perp}\cap y^{\perp}$, then indeed $\A_u\cap\A_v$ and $\A_u\cup\A_v$ are stable under resolution. In the case where one of $u,v\in x^{\perp}\cap y^{\perp}$ and $u,v\neq x,y$, changes under resolution will be controlled by $\overline{\A_u}\setminus\A_u, ~ \overline{\A_v}\setminus\A_v$. This implies that changes in $\overline{\A_u}\cap\overline{\A_v}$ are contained in $(\overline{\A_u}\setminus\A_u)\cap(\overline{\A_v}\setminus\A_v)= (\overline{\A_u}\cap\overline{\A_v})\setminus(\A_u\cup\A_v)$ so, $\A_u\cap\A_v$ and $\A_u\cup\A_v$ are also stable after resolving $xy$. \medskip

Suppose now that $v=x$ and $u\in x^{\perp}\setminus ((x^{\perp}\cap y^{\perp})\cup\{y\})$. Then $\A_u$ and $\overline{\A_u}$ are stable after resolution. From the graph theoretical point of view, it is easy to see that $\overline{S_u}\cap\overline{S_x}$ (as a subgraph of $\Gamma$) remains invariant after resolving the edge $xy$ (considering the same intersection inside $\Gamma_{xy}$). Since $\A_u\cap\A_x\subset \overline{\A_u}\cap\overline{\A_x}$, we deduce that indeed $\A_u\cap\A_x$ and $\A_u\cup\A_x$ are also stable under resolution. The same reasoning holds when $v=y$ and $u\in y^{\perp}\setminus ((x^{\perp}\cap y^{\perp})\cup\{x\})$. This concludes the proof.
\eop

\begin{remark}Notice that the reason {\em why} the previous spaces are not affected by the resolution of the edge $xy$ comes as a direct consequence of the fact that the (loose) subgraphs of $\Gamma$ defining them do not contain the edge $xy$. 
\end{remark}

\begin{remark}
Note that the equality of the last terms in the right-hand sides of (\ref{CLASS}) in Lemma \ref{local} can also be obtained by applying Lemma \ref{affect2}.
\end{remark}

\medskip
\subsection{``Full cones''}

We first handle a useful specific case of graphs.

\begin{lemma}
\label{full}
Suppose that either $y^\perp = (x^\perp \cap y^\perp) \cup \{x\}$, or $x^\perp = (x^\perp \cap y^\perp) \cup \{y\}$. Then $[\mF_k(\Gamma)] \in \Z[\bL]$. 
\end{lemma}
{\em Proof}.\quad
Suppose w.l.o.g. that $y^\perp = (x^\perp \cap y^\perp) \cup \{x\}$; then for all vertices $v$ in $\Gamma$, we have that either $v = x$ or $v \sim x$. 
It follows immediately that
\begin{equation}
[\mF_k(\Gamma)] = [\bA_x] + [\mF_k(x^\perp \cap \Gamma)].
\end{equation}
By induction, applied on the second term in the right-hand side, the lemma follows. \eop \\

\subsection{Without external edges}

We assume that there are no edges $uv$ with $u\in x^{\perp}\setminus ((x^{\perp}\cap y^{\perp})$ \newline $\cup\{y\})$
and $v \in y^{\perp}\setminus ((x^{\perp}\cap y^{\perp})\cup\{x\})$ (call such edges ``external'') | the case where such edges exist will be handled separately below.\\

We also suppose that $y^\perp \ne (x^\perp \cap y^\perp) \cup \{x\}$ and $x^\perp \ne (x^\perp \cap y^\perp) \cup \{y\}$, since otherwise the statement is already true by the previous subsection.

Let $u \ne y$ be any vertex in $x^{\perp} \setminus (x^\perp \cap y^\perp)$; let $e := ux$. 
Let $\Gamma^e$ be the graph $\Gamma$ without the edge $e$ (while not deleting $u$ and $x$); similarly, we define $\Gamma^e_{xy}$. 
As $\Gamma^e$ is a subgraph of $\Gamma$, induction implies that $[\mF_k(\Gamma^e)] \in \Z[\bL]$. Also, by Lemma \ref{under}, $[\mF_k(\Gamma_{xy})] \in
\Z[\bL]$ if and only if $[\mF_k(\underline{\Gamma_{xy}})] \in \Z[\bL]$, and by induction, the latter expression is true since $\underline{\Gamma_{xy}}$ is a 
subgraph of $\Gamma$. In the same way, $[\mF_k(\Gamma^e_{xy})] \in \Z[\bL]$. Now consider $[\mF_k(\Gamma)] - [\mF_k(\Gamma^e)]$. 
Then obviously $[\mF_k(\Gamma)] - [\mF_k(\Gamma^e)] \in \Z[\bL]$ if and only if $[\mF_k(\Gamma)] - [\mF_k(\Gamma_{xu})] \in \Z[\bL]$; by 
Lemma \ref{local}, this holds if and only if $[\mF_k(\Gamma \cap \fP_{x,u})] - [\mF_k(\Gamma_{xu} \cap \fP_{x,u})] \in \Z[\bL]$. Now by our 
assumption, we have 
\begin{equation}
\begin{cases}
\Gamma \cap \fP_{x,u} &\ne\ \Gamma;\\
\Gamma_{xu} \cap \fP_{x,u} &\ne\ \Gamma_{xu},
\end{cases}
\end{equation}
so that induction yields that $[\mF_k(\Gamma)] - [\mF_k(\Gamma^e)] \in \Z[\bL]$. Since $[\mF_k(\Gamma^e)] \in \Z[\bL]$, it follows that $[\mF_k(\Gamma)] \in \Z[\bL]$. \eop \\

\subsection{With external edges}
\label{frizzle2}

Now suppose $\Gamma$ {\em has} external edges. Suppose $\Gamma'$ is the subgraph of $\Gamma$ which one obtains by deleting one chosen
external edge $e = uv$. By induction we know that $[\mF_k({\Gamma}')]$ is in $\Z[\bL]$. 
Then
\begin{equation}
\mF_k(\Gamma) = \mF_k(\Gamma') \coprod \Big((\bA_u \cup \bA_v) \setminus (\cup_{s \in \Gamma'}\bA_{s}')\Big),
\end{equation}
where $\bA_w$ is the local affine space at $w$ in $\mF_k(\Gamma)$, and $\bA_t'$ is the local affine space at $t$ in $\mF_k({\Gamma}')$ (note that $\cup_{s \in \Gamma'}\bA_s' = \mF_k(\Gamma')$). Then
\begin{equation}
[\mF_k(\Gamma)] = [\mF_k(\Gamma')] + \underbrace{\Big[(\bA_u \cup \bA_v) \setminus (\cup_{s \in \Gamma'}\bA_{s}')\Big]}_{(\mathrm{A})}.
\end{equation}

Doing the same for $\Gamma_{xy}$, we obtain that 
\begin{equation}
[\mF_k(\Gamma_{xy})] = [\mF_k(\Gamma'_{xy})] + \underbrace{\Big[(\bA_u \cup \bA_v) \setminus (\cup_{s \in \Gamma'_{xy}}\bA_{s}')\Big]}_{(\mathrm{B})},
\end{equation}
where all the local affine spaces are now considered in $\Gamma_{xy}$ or $\Gamma'_{xy}$.

By Lemma \ref{affect2}, we have that (A) $=$ (B). For, $(\A_u \cup \A_v) \cap \A_x' = \A_u \cap \A_x'$ and $(\A_u \cup \A_v) \cap \A_y'$ $= 
\A_v \cap \A_y'$ do not change when resolving $xy$, and if $w \in x^{\perp} \cap y^{\perp}$, then $(\A_u \cup \A_v) \cap \A_w'$ also does not change through resolution. All the other cases are covered by Lemma \ref{affect2}.



After applying induction, we now get that $[\mF_k(\Gamma)] - [\mF_k(\Gamma')] \in \Z[\bL]$. \eop \\






\medskip
\subsection{End of the proof of Theorem \ref{main2}}

Starting from a connected loose graph $\Gamma$, we first note that if $\Gamma$ is a loose tree, the result follows  from Theorem \ref{D3.1}. So suppose that $\Gamma$ is not a loose tree. Choose any edge $xy$ that is contained in a loose spanning tree $T$, and resolve $xy$. 
In the above we have shown that
\begin{equation}
\label{dec5}
[\mF_k(\Gamma)] - [\mF_k(\Gamma_{xy})] \ \in \ \Z[\bL]. 
\end{equation}

Now there are two ways to proceed.\\

\begin{itemize}
\item[(A)]
Carry out surgery on the scheme $\mF_k(\Gamma)$ (using the loose tree $T$), and eventually wind up with a scheme $\mF_k(\overline{T})$, where $\overline{T}$ is a loose tree containing $T$, cf. \S\S \ref{cplg}. We have seen that $[\mF_k(\overline{T})] \in \Z[\bL]$ in Theorem \ref{D3.1}. Since by (\ref{dec5}) each difference between Grothendieck classes of consecutive steps is an element of $\Z[\bL]$, we conclude that the same is true for the initial class $[\mF_k(\Gamma)]$ as well.
\item[(B)]
Use the induction hypothesis to conclude that $[\mF_k(\Gamma_{xy})] \in \Z[\bL]$, so that $\mF_k(\Gamma) \in \Z[\bL]$.\\
\end{itemize}

\medskip
This concludes the proof of Theorem \ref{main}. \eop \\



\bigskip
\begin{thebibliography}{99}



\bibitem{Deitmarschemes2}
{A.~Deitmar}. 
Schemes over $\mathbb{F}_1$, in {\em Number Fields and Function Fields -- Two Parallel Worlds}, pp. 87--100, Progr. Math. {\bf 239}, Birkh\"{a}user Boston, Boston, MA, 2005. 



\bibitem{Deitmarcongruence}
{A.~Deitmar}. 
Congruence Schemes, {\em International Journal of Mathematics} {\bf 24} (2013), 46pp. 



%
%



%
%



\bibitem{Katz}
T.~Hausel and F.~Rodriguez-Villegas, Fernando.
Mixed Hodge polynomials of character varieties,
With an appendix by Nicholas M. Katz,
{\em Invent. Math.} {\bf 174} (2008), 555--624. 

\bibitem{Kurozeta}
{N.~Kurokawa}. 
 Zeta functions over $\mathbf{F}_1$, {\em Proc. Japan Acad. Ser. A Math. Sci.} {\bf 81} (2005), 180--184.
 
 
 

 
 
 \bibitem{Manin}
{Yu.~Manin}.
Lectures on zeta functions and motives (according to Deninger and Kurokawa), Columbia University Number Theory Seminar (New York, 1992), {\em Ast\'{e}risque} {\bf 228} (1995), 121--163.

 
 \bibitem{MMKT}
 M.~M\'{e}rida-Angulo and K.~Thas. Deitmar schemes, graphs and zeta functions, Revised, 50pp.
 
 \bibitem{MMKT1}
 M.~M\'{e}rida-Angulo and K.~Thas. Automorphisms of Deitmar schemes, I. Functoriality and trees, Submitted, 28pp.  
 
 \bibitem{MMKT2}
  M.~M\'{e}rida-Angulo and K.~Thas.
 The structure of Deitmar schemes, II. Zeta functions and automorphism groups, Submitted, 8pp.
 
%
\bibitem{KT-Japan}
K.~Thas. The structure of Deitmar schemes, I, {\em Proc. Japan Acad. Ser. A Math. Sci.} {\bf 90} (2014),  21--26.

\bibitem{Chap1}
K.~Thas. The Weyl functor | Introduction to Absolute Arithmetic, in {\em Absolute Arithmetic and $\Fun$-Geometry}, Koen Thas (ed.), pp. 3--38, EMS Publishing House, Z\"{u}rich, 2016.
%



\bibitem{anal}
{J.~Tits}.
Sur les analogues alg\'{e}briques des groupes semi-simples complexes, in
{\em Colloque d'Alg\`{e}bre Sup\'{e}rieure, Bruxelles du 19 au 22 d\'{e}c. 1956}, pp. 261-289, \'Etablissements Ceuterick, Louvain; Librairie Gauthier-Villars, Paris,  1957.

\end{thebibliography}
\end{document}